\documentclass{elsarticle}
\usepackage[utf8]{inputenc}

\usepackage{float}
\usepackage[dvipsnames,
svgnames,
cmyk]{xcolor}     
\usepackage{graphicx}         
\usepackage[outercaption]{sidecap}
\usepackage{tabularx}
\usepackage{multirow}       
\usepackage{booktabs}
\usepackage{url}
\usepackage{natbib}
\usepackage{hyperref}

\usepackage{amssymb}    
\usepackage{amsthm}     
\usepackage{thmtools}   
\usepackage{mathtools}  
\usepackage{mathrsfs}   
\usepackage{dsfont}     
\newcommand{\R}{\mathbb{R}}    
\newcommand{\N}{\mathbb{N}}    
\newcommand{\x}{\boldsymbol{x}}
\newcommand{\y}{\boldsymbol{y}}
\newcommand{\z}{\boldsymbol{z}}

\newtheorem{definition}{Definition}
\newtheorem{remark}{Remark}

\newcommand{\bs}{\boldsymbol}

\begin{document}
	
	\begin{frontmatter}
		
		
		
		\title{Variably Scaled Persistence Kernels (VSPKs) for persistent homology applications}
		
		
		\author{S. De Marchi$^{*}$}
		\ead{demarchi@math.unipd.it}
		
		\author{F. Lot$^{**}$}
		\ead{lotf@mathematik.uni-marburg.de}
		
		\author{F. Marchetti$^{*}$}
		\ead{francesco.marchetti@math.unipd.it}
		
		\author{D. Poggiali$^{***}$}
		\ead{dr.davide.poggiali@gmail.com}
		
		\address{$^{*}$Dipartimento di Matematica \lq\lq Tullio Levi-Civita\rq\rq, Universit\`a di Padova, 35121 Padua, Italy;\\$^{**}$Fachbereich Mathematik und Informatik, Philipps-Universit\"at Marburg, 35032 Marburg, Germany;\\ $^{***}$FAR Networks srl, 20063 Cernusco S.N., Milan, Italy}

		\begin{abstract}
			
		In recent years, various kernels have been proposed in the context of \textit{persistent homology} to deal with \textit{persistence diagrams} in supervised learning approaches. In this paper, we consider the idea of variably scaled kernels, for approximating  functions and data, and we interpret it in the framework of persistent homology. We call them {\it Variably Scaled Persistence Kernels (VSPKs)}. These new kernels are then tested in different classification experiments. The obtained results show that they can improve the performance and the efficiency of existing standard kernels.

		\end{abstract}
		
		
	\end{frontmatter}

\section{Introduction}
\label{sec:1}

Let $\Omega \subset \R^v$ and let $\mathcal{X} = \{\x_1,\dots, \x_n \} \subset \Omega$ be a set of input data, $v, n \in \N$. Assume that each element $\x_i$ is uniquely associated to a \textit{label} (or \textit{class}) $y_i \in \mathcal{T}$, with $\mathcal{T} = \{c_1,\dots,c_t\}$, so that we can define the dataset $\mathcal{Z}=\{(\x_i, y_i)\:|\:\x_i\in\mathcal{X},\:y_i\in\mathcal{T}\}$.
The supervised learning task consists in finding a decision function $s: \Omega \longrightarrow \mathcal{T}$ such that:
\begin{enumerate}
	\item[(i)]
	it models the input-output relation in $\mathcal{Z}$;
	\item[(ii)]
	it models the input-output relation of unseen labeled instances\\ $\{(\bs{\xi}_i,y_i)\}_{i=1,\dots,m}\subset(\Omega\setminus \mathcal{X})\times \mathcal{T}$, $m\ge 1$.
\end{enumerate}
The \textit{generalization} capability required in (ii) is fundamental in the learning problem, since it is trivial to find a decision function that satisfies (i); for an introduction concerning (statistical) learning theory, refer to e.g. \cite[\S 1.3]{scholkopf2002learning} and \cite{vapnik1999overview}.

Kernel methods are well-established tools in supervised machine learning, as well as in a variety of research and applied fields \cite{schaback_wendland_2006,shawe2004kernel}. The flexibility provided by kernel-based schemes allows the handling of different kinds of possible structured data, e.g, graphs and words, which are encoded in some dot product space where even complex patterns may be distinguished \cite{borgwardt2005protein,collins2001convolution,lampert2009kernel,vishwanathan2010graph}.

In the following, our data consist of \textit{persistence diagrams}, which represent an output of the so-called \textit{persistent homology}, in the framework of Topological Data Analysis (TDA). The usage of topological methods and analysis techniques to extract significant features and patterns from data is receiving more and more interest and, in particular, persistent homology captures the evolution of topological features in the data; for a complete overview concerning TDA, see e.g. \cite{carlsson2009topology,Edelsbrunner2000TopologicalPA, fomenko2012visual}. In the last years, several kernels specifically devoted to deal with the peculiar structure of persistence diagrams have been proposed \cite{carriere2017sliced,kusano2017kernel,le2018persistence,reininghaus2015stable}, therefore the construction of suitable kernels is a very active research line in this context.

Here, we introduce what we call the Variably Scaled Persistence Kernels (VSPKs), which are inspired by the variably scaled kernels introduced in the context of approximation theory in \cite{bozzini2015interpolation} and employed as a feature augmentation strategy in \cite{learningvsk,vskieee} for kernel-based learning. Indeed, our aim is to build a bridge between classical variably scaled kernels and kernels defined in the context of persistent homology.
After defining VSPKs, we design some scaling functions and we test the performance of the resulting kernels by experimenting with three different datasets. The obtained results show that the variably scaled setting may yields to better classification outcomes and efficiency with respect to the standard setting, and thus can be considered for further investigations and applications involving persistence diagrams.

The paper is organised as follows. In Sections \ref{sec:2} and \ref{sec:3}, we introduce kernel-based learning with Support Vector Machines (SVMs) and we recall some notions concerning persistent homology, also presenting some kernels of recent introduction that are dedicated to persistent diagrams. VSPKs are proposed and analysed in Section \ref{sec:4}, and related numerical experiments are exposed in Section \ref{sec:5}. Finally, in Section \ref{sec:6} we present some conclusions and final remarks.

\section{Positive definite and variably scaled kernels}
\label{sec:2}

Let $\kappa: \Omega \times \Omega  \longrightarrow \R$ be a kernel. Given a set of data $\mathcal{X}=\{\x_1,\dots,\x_n\} \subset \Omega$, the $n \times n$ matrix $K$ with elements $K_{i j} := \kappa(\x_i, \x_j)$, $i,j=1,\dots,n$, is the \textit{Gram matrix} of the kernel $\kappa$ with respect to $\mathcal{X}$. If $\kappa$ is positive definite (strictly positive definite) on $\Omega\times\Omega$, i.e., $K$ is positive semi-definite (definite) for all possible datasets in $\Omega$, then it is possible to decompose the kernel according to Mercer theorem \cite{mercer}, and to interpret such decomposition as an inner product in a \textit{Reproducing Kernel Hilbert Space} (RKHS) $\mathscr{F}$. Indeed, there exists a (non-unique) \textit{feature map} $\Phi: \Omega\longrightarrow \mathscr{F}$ such that
\begin{equation*}
	\kappa(\x, \y) = \langle \Phi(\x), \Phi(\y) \rangle_\mathscr{F} \qquad  \x, \y \in \Omega,
\end{equation*}
being $\langle \cdot, \cdot \rangle_\mathscr{F}$ the bilinear form related to the RKHS (e.g. $\Phi(\x)=\kappa(\cdot,\x)$). Moreover, the kernel $\kappa$ induces a distance $d_\kappa$ on $\Omega$
\begin{equation}\label{eq:induced}
	d_\kappa(\x, \y) := \kappa(\x, \x) + \kappa(\y, \y) - 2\kappa(\x, \y).
\end{equation}

Variably Scaled Kernels (VSKs) have been introduced in \cite{bozzini2015interpolation} in the context of kernel-based approximation, with the aim of overcoming instability issues. Then, they have been extended to work in a more general setting in \cite{learningvsk}, as presented in the following form. Let $\Lambda \subseteq \mathbb{R}^{\nu}$, $\nu>0\in \mathbb{N}$ and let $\kappa: \tilde{\Omega}\times \tilde{\Omega}\longrightarrow\R$ be
a continuous (strictly) positive definite kernel, where $\tilde{\Omega}=\Omega \times \Lambda \subseteq \mathbb{R}^{v+\nu}$. Given a \textit{scaling} function $\Psi: \Omega \longrightarrow \Lambda,$
a VSK $\kappa_{\Psi}: \Omega  \times \Omega \longrightarrow \mathbb{R}$ is defined as
\begin{equation}\label{defvskn}
	\kappa_{\Psi}(\bs{{x}},\bs{{y}})= \kappa((\bs{{x}},\Psi(\bs{x})),(\bs{y},\Psi(\bs{y})))
\end{equation}
for $\bs{x},\bs{y}\in\Omega$. The function $\Psi$ can be interpreted as a \textit{feature augmentation} map, which adds $\nu$ coordinates (features) to the original sample. In this view, the VSK setting has been analysed in \cite{learningvsk} as a \textit{stacking} technique, which is capable of enhancing the prediction performances of classical kernel-based classifiers such as, e.g., SVMs.

Letting $\boldsymbol{x}=(x_1,\dots,x_d)^{\intercal} \in \Omega$, we recall that a binary (i.e. $\mathcal{T}=\{-1,+1\}$) SVMs classifier is characterised by the decision function
\begin{equation*}
	s(\boldsymbol{x})= \textrm{sign}( h(\boldsymbol{x}))= \textrm{sign}( \langle \Phi(\bs{x}),\boldsymbol{w}\rangle_{\mathscr{F}}+b),
\end{equation*}
where
\begin{equation*}\label{ww}
	\boldsymbol{w}=  \sum_{i=1}^n   \alpha_i y_i \Phi(\boldsymbol{x}_i)\in {\mathscr{F}}.
\end{equation*}
The coefficients $\bs{\alpha}=(\alpha_1,\dots,\alpha_n)\in\R^n$ are the solution of the following \textit{soft margin} problem \cite[\S 18, p. 346--347]{Fasshauer15}
\begin{equation*} 
	\left\{
	\begin{array}{l}
		\min_{\boldsymbol{\alpha} \in \mathbb{R}^{n} } \dfrac{1}{2} \sum_{i=1}^{n} \sum_{j=1}^{n} \alpha_i \alpha_j y_i y_j \kappa(\boldsymbol{x}_i,\boldsymbol{x}_j)- \sum_{i=1}^{n} \alpha_i,  \\
		\textrm{s.t.  }   \sum_{i=1}^{n} \alpha_i y_i =0,\\
		0 \leq \alpha_i \leq \zeta, \hskip 0.3cm i=1,\ldots,n,
	\end{array}
	\right.
\end{equation*}
where $[0,\zeta]^n$ is the \textit{bounding box}, with $\zeta \in  [0,+\infty)$. Usually, a binary SVMs classifier is extended to the multiclass setting by considering a \textit{one-vs-rest} approach.

\section{Persistent homology and kernels}
\label{sec:3}

\subsection{Basics on persistent homology}
In the following, we recall some basic ideas about some tools of persistent homology. For a more detailed treatment, especially concerning the algebraic aspects of the construction, we refer e.g. to \cite{Edelsbrunner2000TopologicalPA,Guillemard2017}.

Let our data $\mathcal{X}=\{\x_1,\dots,\x_n\} \subset \Omega$ be interpreted as a set of vertices sampled from some manifold $M$, and suppose that we wish to highlight some of its intrinsic homological properties. Letting $\varepsilon>0$, a possible concrete way to proceed consists of studying
$$
M^{\mathcal{X},\varepsilon} = \bigcup_{i=1}^n B(\bs{x}_i,\varepsilon)
$$
as an approximation of $M$, where $B(\bs{x}_i,\varepsilon)$ is the ball of radius $\varepsilon$ and centre $\bs{x}_i$. We can associate to $M^{\mathcal{X},\varepsilon}$ a \textit{Vietoris-Rips simplicial complex} $K^{\mathcal{X},\varepsilon}$, where two distinct vertices $\bs{x}_i,\bs{x}_j$ are connected by an edge if and only if $\lVert \bs{x}_i-\bs{x}_j\lVert_2\le\varepsilon$. Formal linear combinations of $r$-dimensional faces in $K^{\mathcal{X},\varepsilon}$ form the \textit{$r$-chains} group $C_r^{\mathcal{X},\varepsilon}$. Moreover, letting $[\bs{x}_{i_0},\dots,\bs{x}_{i_{r}}]$ be the face constructed upon distinct vertices $\bs{x}_{i_0},\dots,\bs{x}_{i_{r}}$ in $\mathcal{X}$, $i_0,\dots,i_r\in\{1,\dots,n\}$, we define the linear \textit{boundary operator} $\partial_r:C_r^{\mathcal{X},\varepsilon}\longrightarrow C_{r-1}^{\mathcal{X},\varepsilon}$ as
$$
\partial_r[\bs{x}_{i_0},\dots,\bs{x}_{i_{r}}]=\sum_{j=0}^r(-1)^r[\bs{x}_{i_0},\dots,\bs{x}_{i_{j-1}},\bs{x}_{i_{j+1}},\dots,\bs{x}_{i_{r}}].
$$
The \textit{$r$-cycles} and \textit{$r$-boundaries} groups are defined as $Z_r^{\mathcal{X},\varepsilon}=\mathrm{ker}\partial_r$ and $B_r^{\mathcal{X},\varepsilon}=\mathrm{im}\partial_{r+1}$, respectively. Furthermore, the rank of the \textit{$r$-homology} group $H_r^{\mathcal{X},\varepsilon}=Z_r^{\mathcal{X},\varepsilon}/B_r^{\mathcal{X},\varepsilon}$ expresses the concept of \textit{$r$-dimensional holes} in $K^{\mathcal{X},\varepsilon}$.

Finding an optimal $\varepsilon^{\star}$ that represents the intrinsic geometric properties of $M$ is a tough and unstable process. Instead, one may analyse the whole \textit{filtration} $\{M^{\mathcal{X},\varepsilon}\:|\:\varepsilon>0\}$. In particular, letting $\varepsilon_1<\dots<\varepsilon_u$ be increasing real numbers, the nested sequence $K^{\mathcal{X},\varepsilon_1}\subseteq\dots\subseteq K^{\mathcal{X},\varepsilon_u}$ is obtained. Then, for $r\ge0$ and $i\in\{1,\dots,u\}$, we consider the \textit{$\ell$-persistent homology} group 
$$
H_{r,\ell}^{\mathcal{X},\varepsilon_i}=Z_r^{\mathcal{X},\varepsilon_i}/(Z_r^{\mathcal{X},\varepsilon_i}\cap B_r^{\mathcal{X},\varepsilon_{i+\ell}}).
$$
The group $H_{r,\ell}^{\mathcal{X},\varepsilon_i}$ contains the homology classes that persist in the \textit{time interval} $[i,i+\ell]$, i.e., they are \textit{born} with $K^{\mathcal{X},\varepsilon_j}$ for some $j<i$ and they are \textit{alive} with $K^{\mathcal{X},\varepsilon_{i+\ell}}$. We point out that such homology classes might persist \textit{indefinitely} (we will denote this case as an $\infty$ level), or they might \textit{die} with a certain $\varepsilon_j$, $i<j\le u$.

\begin{figure}[h]
	\centering
	\includegraphics[width=.49\linewidth]{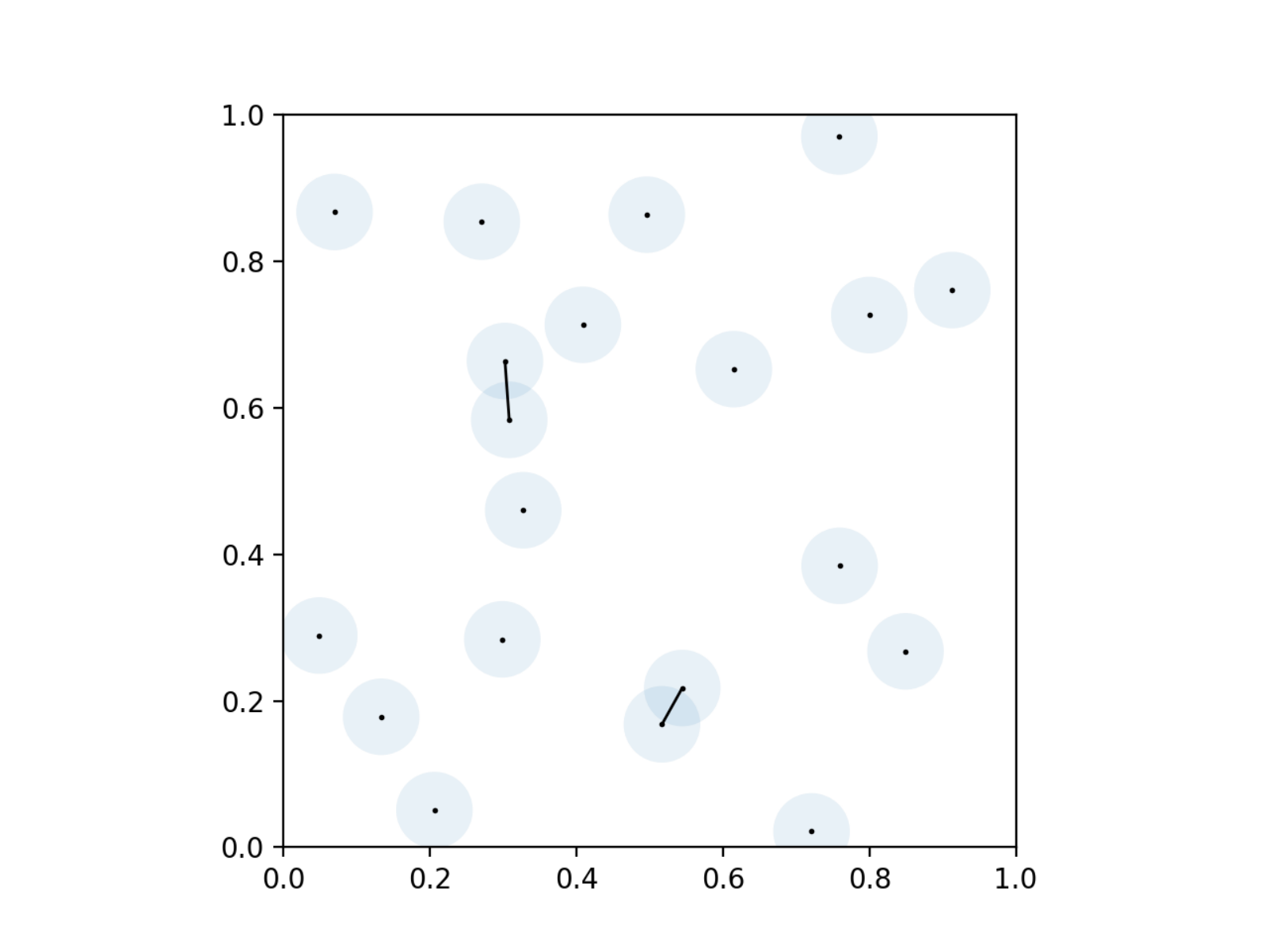}
	\includegraphics[width=.49\linewidth]{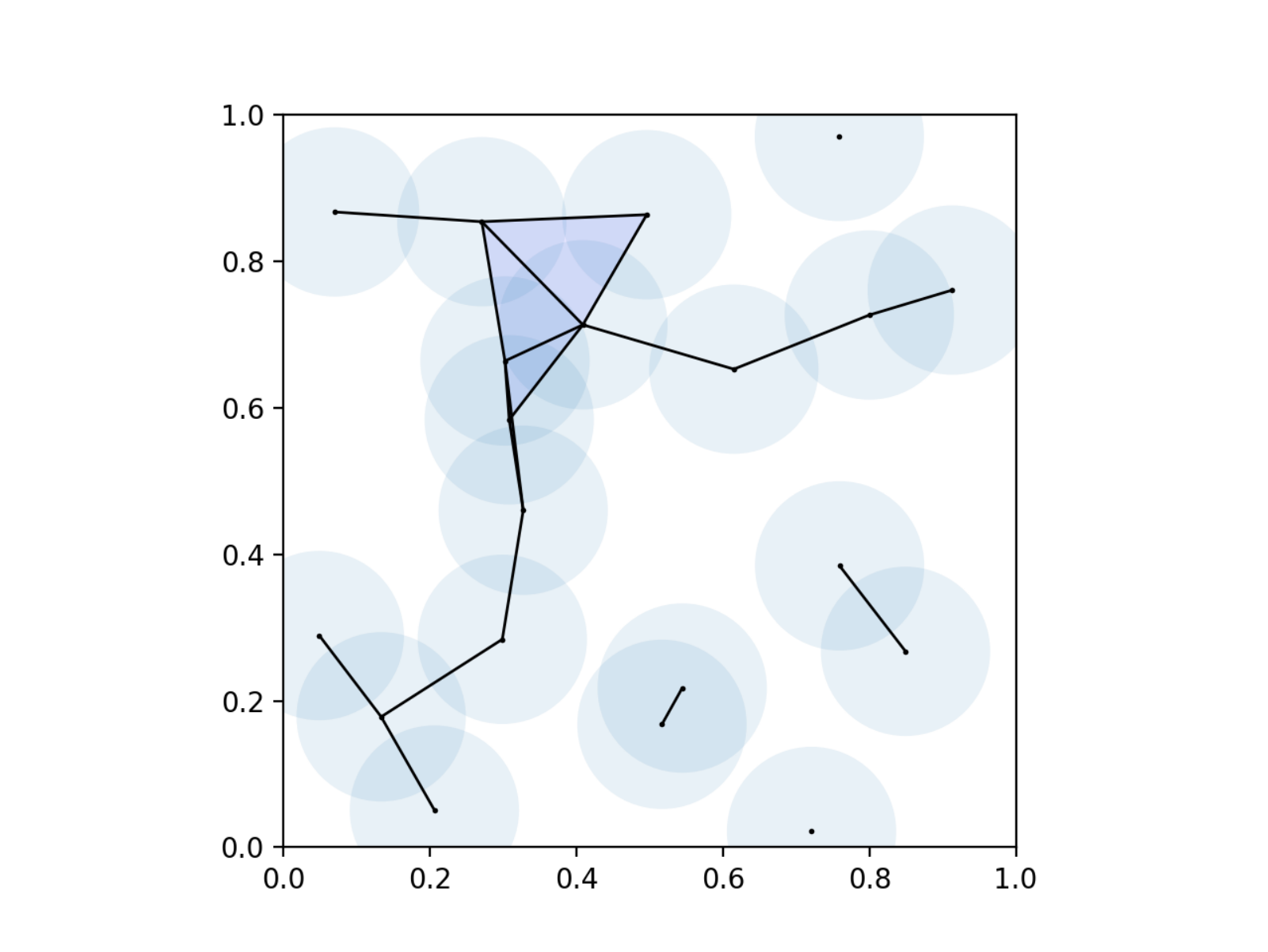} \\
	\includegraphics[width=.49\linewidth]{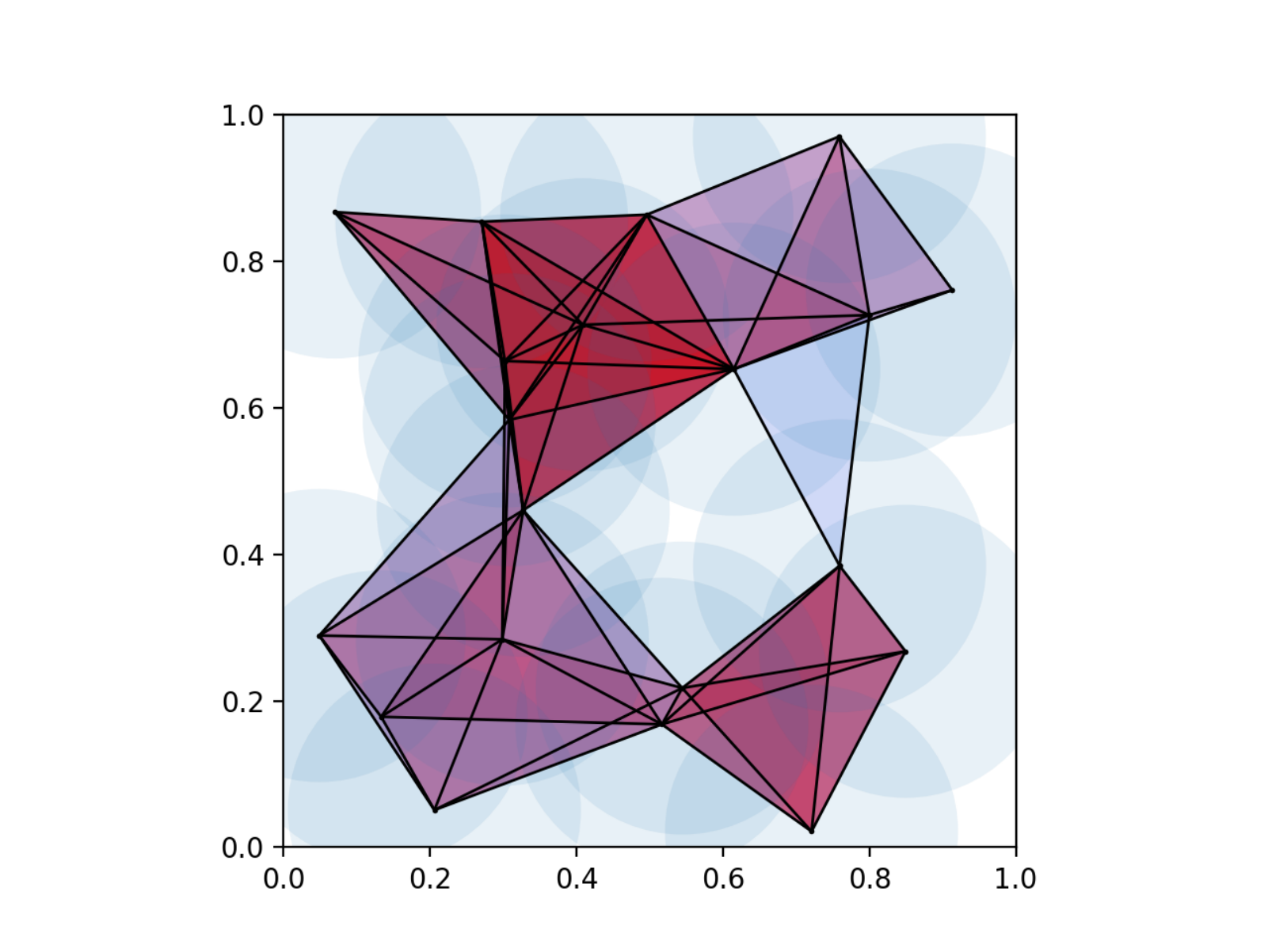}
	\includegraphics[width=.49\linewidth]{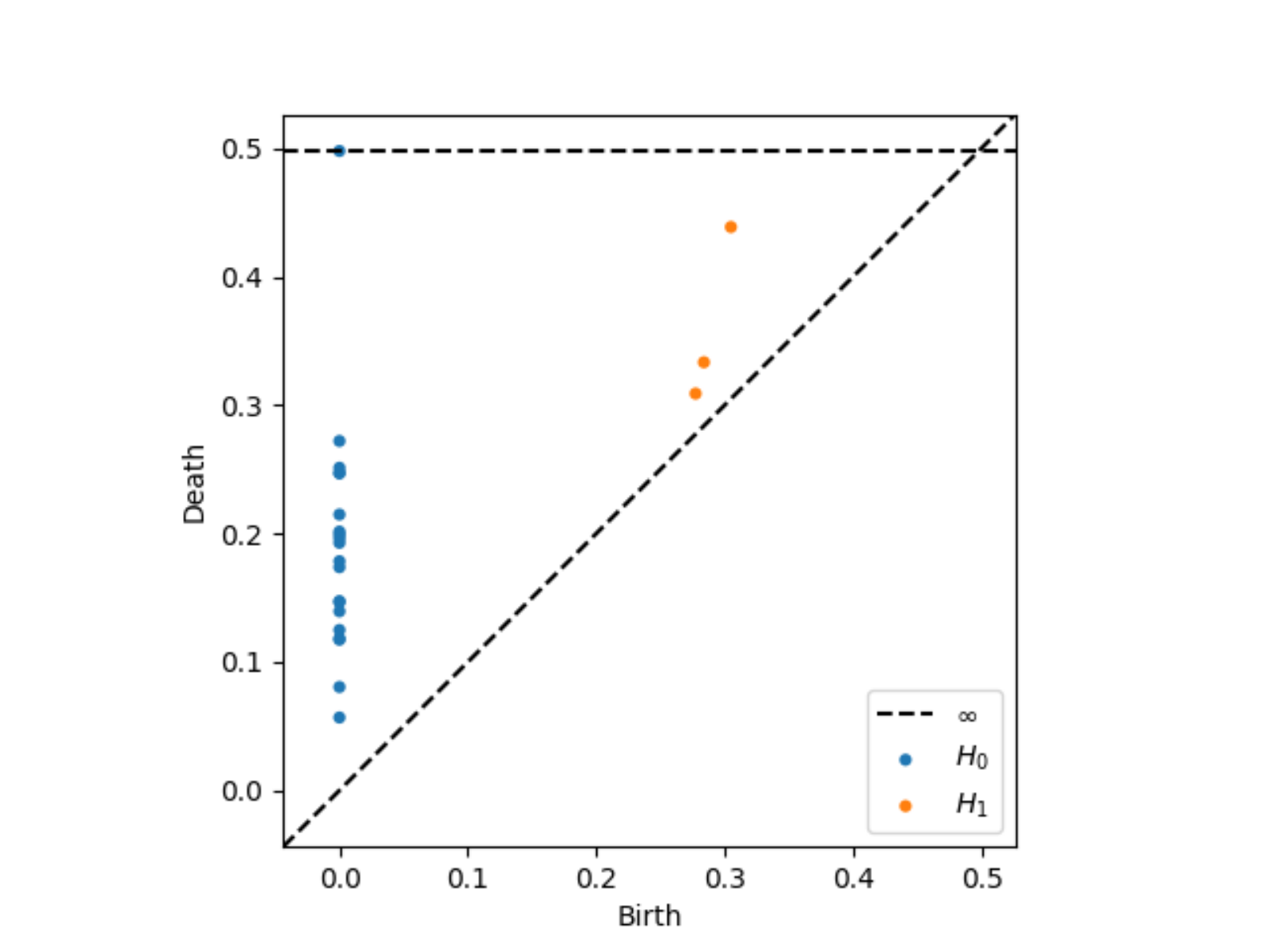}
	\caption{Given a random set of points $\mathcal{X}$, the construction of the Vietoris-Rips complex along the filtration is depicted on top and bottom left figures. The persistence diagrams $D_0(\mathcal{X},\bs{\varepsilon})$ and $D_1(\mathcal{X},\bs{\varepsilon})$ are overlapped in the bottom right figure.}
	\label{fig:example_pd_point}
\end{figure}

Therefore, each element of the persistent homology groups obtained by considering the whole filtration can be represented by a birth-death pair $(b,d)\in\R^2_{+}$, $b=\varepsilon_{h},\;d=\varepsilon_{k}$ for some $h\in\{1,\dots,m\}$, $k\in\{1,\dots,m\}\cup\{\infty\}$, $h<k$. We say that $d-b$ is the \textit{persistence} of $(b,d)$, and that a birth-death pair is $r$-dimensional if it is related to $r$-dimensional homology groups. Moreover, letting $\bs{\varepsilon}=(\varepsilon_1,\dots,\varepsilon_u)$, we define a \textit{persistence diagram} $D_r(\mathcal{X},\bs{\varepsilon})$ related to the filtration $K^{\mathcal{X},\varepsilon_1}\subseteq\dots\subseteq K^{\mathcal{X},\varepsilon_m}$ as
\begin{equation}\label{eq:pdiagram}
	D_r(\mathcal{X},\bs{\varepsilon})=\{(b,d)\:|\:(b,d)\in P_r(\mathcal{X},\bs{\varepsilon})\}\cup \mathcal{B},\quad r\ge 0,
\end{equation}
where $P_r(\mathcal{X},\bs{\varepsilon})$ denotes the set of $r$-dimensional birth-death pairs obtained with the filtration and $\mathcal{B}=\{(z,z)\:|\: z\ge 0\}$.
We remark that $D_r(\mathcal{X},\bs{\varepsilon})$ is a multiset, since a couple $(b,d)$ might appear more than once, i.e., might have multiplicity greater than one. Furthermore, the bisector $\mathcal{B}$ is composed by an infinite number of elements characterised by infinite multiplicity, and it is added in order to achieve some \textit{uniformity} among different persistence diagrams and facilitate the formulation of proper metrics, as we present below. In Figure \ref{fig:example_pd_point}, we display the steps and the properties of the discussed analysis by means of a two-dimensional example. We observe that such analysis captures some intrinsic geometrical properties of $\mathcal{X}$. In particular, the unique most persistent $0$-dimensional pair, which is significantly far from the others, suggests that a unique connected component underlies $\mathcal{X}$. Moreover, the persistent $1$-dimensional pair indicates the presence of a $1$-dimensional hole in the structure of the dataset, which is highlighted by the simplicial complex depicted in the bottom left figure.

Persistence diagrams show some stability properties with respect to perturbations of the involved dataset \cite{cohen2007stability}. To better clarify this aspect, let us recall some useful metrics. Let $\mathcal{X}, \mathcal{Y} \subset \Omega$ be two non-empty datasets. The \textit{Hausdorff distance} is defined as
\begin{equation*}
	d_H(\mathcal{X},\mathcal{Y}) = \max \bigg\{\sup_{x \in \mathcal{X}} \inf_{y \in \mathcal{Y}} \lVert x-y\lVert_{\infty}, \sup_{y \in \mathcal{Y}} \inf_{x \in \mathcal{X}} \lVert x-y\lVert_{\infty}\bigg\}.
\end{equation*}
Moreover, letting $D_r(\mathcal{X},\bs{\varepsilon})$ and $D_r(\mathcal{Y},\bs{\varepsilon})$ be persistence diagrams for some $r\ge 0$ and filtration radii vector $\bs{\varepsilon}$, we recall the \textit{$p$-Wasserstein} distance
\begin{equation*}
	d_{W,p}(D_r(\mathcal{X},\bs{\varepsilon}),D_r(\mathcal{Y},\bs{\varepsilon})) = \bigg(\inf_{\gamma\in\Gamma} \sum_{x \in D_r(\mathcal{X},\bs{\varepsilon})} \lVert x-\gamma(x)\lVert^p_{\infty}\bigg)^{\frac{1}{p}},
\end{equation*}
where $\Gamma=\{\gamma: D_r(\mathcal{X},\bs{\varepsilon})\longrightarrow D_r(\mathcal{Y},\bs{\varepsilon})\:|\: \gamma\textrm{ is a bijection}\}$. In particular, letting $p\rightarrow \infty$, we obtain the \textit{bottleneck distance}
\begin{equation*}
	d_{W,\infty}(D_r(\mathcal{X},\bs{\varepsilon}),D_r(\mathcal{Y},\bs{\varepsilon}))=d_{B}(D_r(\mathcal{X},\bs{\varepsilon}),D_r(\mathcal{Y},\bs{\varepsilon})) = \inf_{\gamma\in\Gamma} \sup_{x \in D_r(\mathcal{X},\bs{\varepsilon})} \lVert x- \gamma(x)\lVert_{\infty}.
\end{equation*}
We have the stability result \cite{kusano2017kernel}
$$
d_{B}(D_r(\mathcal{X},\bs{\varepsilon}),D_r(\mathcal{Y},\bs{\varepsilon})) \le d_H(\mathcal{X},\mathcal{Y}),
$$
i.e., the bottleneck distance between persistence diagrams is controlled as long as the underlying datasets are close in the Hausdorff metric. 

\subsection{Kernels for persistence diagrams}

In order to better measure similarities between persistence diagrams, various positive definite kernels that are suitable for dealing with the peculiar structure of persistence diagrams have been introduced and studied in the recent literature. In our applications, we will consider the following ones.

\begin{itemize}
	\item 
	The \textbf{Persistence Scale Space (PSS) kernel} \cite{reininghaus2015stable} 
	\begin{equation*}
		\label{eq:pss}
		\kappa_\sigma(D_1,D_2) = \frac{1}{8 \pi \sigma} \sum_{\substack{\y \in D_1 \\ \z \in D_2}} \exp\bigg(-\frac{||\y-\z||^2}{8 \sigma}\bigg) - \exp\bigg(-\frac{||\y-\bar{\z}||^2}{8 \sigma}\bigg),
	\end{equation*}
	
	which is 1-Wasserstein stable, i.e.:
	\begin{equation*}
		\|\kappa_\sigma(D_1,D_2)\|_{L_2(\Omega)} \le \frac{1}{2 \sqrt{\pi} \sigma} d_{W,1}(D_1,D_2)
	\end{equation*}
	
	\item
	The \textbf{Persistence Weighted Gaussian (PWG) kernel} \cite{kusano2017kernel}
	
	\begin{equation*}
		\label{eq:pwg}
		\kappa_G(D_1,D_2;\kappa,\omega) := \exp\bigg(- \frac{1}{2 \tau^2} \| E_\kappa(\mu_{D_1}^{\omega}) - E_\kappa(\mu_{D_2}^{\omega})\|^2_{\mathscr{F}}\bigg) \quad \tau > 0,
	\end{equation*}
	
	which is built upon a \textit{standard} kernel $\kappa$ and a weight function $\omega$, where 
	
	\begin{equation*}
		E_\kappa(\mu_{D_1}^{\omega}) := \sum_{\x \in D_1} \omega(\x) \, \kappa(\cdot, \x).
	\end{equation*}
	
	It is both 1-Wasserstein and bottleneck stable, if we choose as weight function
	\begin{equation*}
		\omega_{\text{arc}}(\x) = \arctan(C (d-b)^{\delta}) \quad \x=(b,d),\, C > 0, \, \delta \in \mathbb{Z}_{>0}.  
	\end{equation*}
	
	Indeed, there exist $\delta\in \mathbb{Z}_{>0}$ and $L>0$ such that
	
	\begin{equation*}
		\| E_\kappa(\mu_{D_1}^{\omega_{\text{arc}}}) - E_\kappa(\mu_{D_2}^{\omega_{\text{arc}}}) \| \le L  d_{B}(D_1,D_2).
	\end{equation*}
	
	\item
	
	The \textbf{Sliced Wasserstein (SW) kernel} \cite{carriere2017sliced}
	
	\begin{equation*}
		\kappa_{SW}(D_1,D_2) := \exp\left(-\frac{SW(D_1,D_2)}{2 \sigma^2}\right).
	\end{equation*}
	
	which is based on the so called Sliced Wasserstein distance, which is equivalent to the 1-Wasserstein distance, i.e.,
	
	\begin{equation*}
		\frac{1}{2M} d_{W,1}(D_1,D_2) \le d_{SW}(D_1,D_2) \le 2 \sqrt{2} \, d_{W,1}(D_1,D_2)
		\label{eq:SWequal}
	\end{equation*}
	
	for some positive constant $M$.
	
\end{itemize}
\begin{remark}\label{rem:specials}
	We observe that the PWG and SW kernels are in the form
	\begin{equation*}
		\kappa(D_r(\mathcal{X},\bs{\varepsilon}),D_r(\mathcal{Y},\bs{\varepsilon})) = \exp \left(- \beta d(D_r(\mathcal{X},\bs{\varepsilon}),D_r(\mathcal{Y},\bs{\varepsilon}) \right))
		\label{eq:kind}
	\end{equation*}
	for some $\beta>0$, where $d(\cdot, \cdot)$ is the distance induced by the underlying metric in the case of the SW (for more details, we refer to the cited seminal papers), while for the PWG kernel the distance is induced by the kernel (see \eqref{eq:induced}).
\end{remark}

\section{Variably Scaled Persistence Kernels}
\label{sec:4}

In the following, our purpose is to interpret the idea underlying
VSKs in the context of persistent homology. The main difference between standard kernels and kernels for persistence diagrams is the structure of the input data. Since persistent diagrams consists of a collection of topological features, i.e. birth-death couples in $\R^2_{\ge\mathcal{B}}=\{\bs{x}=(b,d)\in\R^2_+\:|\: d\ge b\}$, introducing a scaling function whose output lies outside $\R^2_{\ge\mathcal{B}}$ would be meaningless. Hence, letting $\mathfrak{D}_r(\bs{\varepsilon})=\{D_r(\mathcal{X},\bs{\varepsilon})\:|\:\mathcal{X}\subset\Omega\}$, we propose the following definition.

\begin{definition}
	Let $\kappa : \mathfrak{D}_r(\bs{\varepsilon}) \times \mathfrak{D}_r(\bs{\varepsilon}) \longrightarrow \R$ be a kernel for persistence diagrams and let $\Psi: \mathfrak{D}_r(\bs{\varepsilon}) \rightarrow \mathfrak{D}_r(\bs{\varepsilon})$. A variably scaled persistence kernel $\kappa_\Psi$ on $\mathfrak{D}_r(\bs{\varepsilon}) \times \mathfrak{D}_r(\bs{\varepsilon})$ is defined as
	\begin{equation*}
		\kappa_\Psi(D_r(\mathcal{X},\bs{\varepsilon}),D_r(\mathcal{X},\bs{\varepsilon})) := \kappa( \Psi(D_r(\mathcal{X},\bs{\varepsilon})), \Psi(D_r(\mathcal{Y},\bs{\varepsilon})))
	\end{equation*}
	for $D_r(\mathcal{X},\bs{\varepsilon}),D_r(\mathcal{Y},\bs{\varepsilon}) \in \mathfrak{D}_r(\bs{\varepsilon})$.
	\label{def:vspk}
\end{definition}

As in other contexts, a proper \textit{scaling} function $\Psi$ needs to be designed. Let $\widetilde{D_r}(\mathcal{X},\bs{\varepsilon})=D_r(\mathcal{X},\bs{\varepsilon})\setminus \mathcal{B}$ and let $\psi: \mathfrak{D}_r(\bs{\varepsilon}) \rightarrow \R^2_{\ge\mathcal{B}}$ be an auxiliary function defined as
\begin{equation}\label{eq:scaling_fun}
	\psi(D_r(\mathcal{X},\bs{\varepsilon})) = \frac{1}{W}\sum_{\x \in \widetilde{D_r}(\mathcal{X},\bs{\varepsilon})} w(\x) \, \x
\end{equation}
where $w:\R_+^2\longrightarrow \R_+$ is a weight function and $W = \sum_{\x \in \widetilde{D_r}(\mathcal{X},\bs{\varepsilon})} w(\x)$. We observe that $\widetilde{D_r}(\mathcal{X},\bs{\varepsilon})$ contains a finite number of elements (generators), therefore the sum in \eqref{eq:scaling_fun} is always defined (see also Remark \ref{rem:somma_infinita}). We propose the following alternative scaling functions $\Psi$.

\begin{enumerate}
	\item 
	We define $$\Psi_a(D_r(\mathcal{X},\bs{\varepsilon}))=D_r(\mathcal{X},\bs{\varepsilon}) \cup \psi(D_r(\mathcal{X},\bs{\varepsilon})).$$
	\item
	Letting $\rho\in\mathbb{N}$, we first define the set $\widetilde{D_r}(\mathcal{X},\bs{\varepsilon},\rho)$ which consists of the $\rho$ most persistent elements in $\widetilde{D_r}(\mathcal{X},\bs{\varepsilon})$. Then, we define the function $$\Psi_{\rho}(D_r(\mathcal{X},\bs{\varepsilon}))=\widetilde{D_r}(\mathcal{X},\bs{\varepsilon},\rho) \cup \psi\big(D_r(\mathcal{X},\bs{\varepsilon})\setminus \widetilde{D_r}(\mathcal{X},\bs{\varepsilon},\rho)\big)\cup\mathcal{B}.$$
	We remark that $(D_r(\mathcal{X},\bs{\varepsilon})\setminus \widetilde{D_r}(\mathcal{X},\bs{\varepsilon},\rho))\in \mathfrak{D}_r(\bs{\varepsilon})$, therefore $\Psi_{\rho}(D_r(\mathcal{X},\bs{\varepsilon}))$ is well defined.
\end{enumerate}
We observe that $\Psi_a$ plays the role of a \textit{feature augmenting} map, since an additional generator is added in the persistence diagram. On the other hand, $\Psi_{\rho}$ performs a \textit{feature extraction} procedure. Indeed, the resulting persistence diagram consists of the most $\rho$ persistent elements, $\mathcal{B}$ and the remaining elements, which are possible large in numbers, are \textit{compressed} into a single generator.

\begin{remark}
	Since $D_r(\mathcal{X},\bs{\varepsilon})$ contains an infinite number of elements for all $\mathcal{X}\subset\Omega$ by definition, the VSPK $\kappa_{\Psi}$ is still well defined on $\mathfrak{D}_r(\bs{\varepsilon}) \times \mathfrak{D}_r(\bs{\varepsilon})$. Moreover, if $\kappa$ is (strictly) positive definite, so it is $\kappa_{\Psi}$.
\end{remark}

\begin{remark}
	Referring to Remark \ref{rem:specials}, if $\kappa$ is a PWG or SW kernel, then $\kappa_{\Psi}$ can be directly expressed in terms of the distance $d_{\Psi}(\cdot, \cdot)$ induced in the variably scaled setting.
\end{remark}

For the auxiliary function $\psi$, we propose the following weights.
\begin{enumerate}
	\item 
	Let $w_1(\x)=1/|\widetilde{D_r}(\mathcal{X},\bs{\varepsilon})|$, being $|\widetilde{D_r}(\mathcal{X},\bs{\varepsilon})|$ the cardinality of the multiset, i.e., each element is counted with its multiplicity. We denote as \textit{centre of uniform mass} the resulting auxiliary function
	\begin{equation*}
		\label{eq:vspk_coum}
		\psi_1(D_r(\mathcal{X},\bs{\varepsilon})) =  \frac{1}{|\widetilde{D_r}(\mathcal{X},\bs{\varepsilon})|}\sum_{\x \in \widetilde{D_r}(\mathcal{X},\bs{\varepsilon})} \x.
	\end{equation*}
	\item
	Let $\x=(b,d)\in\widetilde{D_r}(\mathcal{X},\bs{\varepsilon})$, where $b,d\in\R_+$ are the birth-death \textit{time} of the element $\x$ (see \eqref{eq:pdiagram}), and let $w_2(x)=d-b$ be the persistence of $\x$. We denote as \textit{centre of persistence} the auxiliary function
	
	\begin{equation*}
		\label{eq:vspk_cop}
		\psi_2(D_r(\mathcal{X},\bs{\varepsilon})) =  \frac{1}{\sum_{\x=(b,d) \in \widetilde{D_r}(\mathcal{X},\bs{\varepsilon})} (d-b)}\sum_{\x=(b,d) \in \widetilde{D_r}(\mathcal{X},\bs{\varepsilon})} (d-b) \, \x.
	\end{equation*}
\end{enumerate}

While the centre of uniform mass is the \textit{barycentre} of the elements of the multiset, the centre of persistence assigns different weights according to the persistence of the elements. This is a natural choice to be analysed, since elements with low persistence are more likely to be related to \textit{noise} structures resulting in the filtration, while elements of large persistence are linked to more representative geometrical features of the dataset (see the example in Figure \ref{fig:example_pd_point}).

\begin{remark}\label{rem:somma_infinita}
	If we take $\x\in D_r(\mathcal{X},\bs{\varepsilon})$ in \eqref{eq:scaling_fun}, additional conditions on the weight function $w$ are needed to guarantee the convergence of the sum. However, such infinite setting is not meaningful to be analysed, since elements in the bisector carry no topological information concerning the dataset. As a further observation, the centre of persistence might be formally computed by summing over $D_r(\mathcal{X},\bs{\varepsilon})$, as in this case $w(\x)=0$ for $\x\in \mathcal{B}$.
\end{remark}

\section{Experiments}
\label{sec:5}

Our aim is to show how the SVM classifier may benefit of the introduced VSPKs.

In the experiments the kernels are handled using Python 3.8 and the modulus scikit-learn \cite{scikit-learn} on a 2.6 Ghz Dual-Core Intel Core i5. Persistence diagrams are constructed via the modulii \textit{persim} \cite{scikittda2019}, \textit{ripser} \cite{ctralie2018ripser} and \textit{GUDHI} \cite{gudhi:urm}. Free and open source PYTHON software is available at
\begin{center}
	\texttt{\url{https://github.com/reevost/vspk_paper}}\:.
\end{center}

We validate the following hyperparameters. About the considered kernels, we follow the guidelines provided by the authors of the seminal papers.
\begin{itemize}
	\item 
	Concerning the SVM classifier, we validate $\zeta\in\{10^{j}\:|\: j=-3,\dots,3\}$.
	\item
	Concerning the PSS kernel, we take $\sigma\in\{10^{j}\:|\: j=-3,\dots,3\}\cup\{5\cdot10^{j}\:|\: j=-3,\dots,2\}$. 
	\item
	Concerning the PWG kernel, the parameters $C$ and $\tau$ of the PWGK are chosen in $\{10^{j}\:|\: j=-2,\dots,2\}$, while $\delta$ is set to $10$ (see \cite[Theorem 3.2]{kusano2016persistence}). Moreover, as underlying standard kernel we use the Gaussian.
	\item
	Concerning the SW kernel, $\sigma$ is obtained following the procedure carried out in \cite[\S 4]{carriere2017sliced}.
\end{itemize}

To assess the performance of the classifiers, we consider the following scores.

\begin{itemize}
	\item 
	The accuracy score
	{\small
		\begin{equation*}
			\text{accuracy} = \frac{\text{true positives}+\text{true negatives}}{\text{true positives}+\text{true negatives}+\text{false positives}+\text{false negatives}}.
		\end{equation*}
		\item
		The f$_1$-score
		\begin{equation*}
			\text{f}_1\text{-score} = 2 \cdot \frac{\text{precision} \cdot \text{recall}}{\text{precision} + \text{recall}},
		\end{equation*}
		where
		\begin{equation*}
			\text{precision} = \frac{\text{true positives}}{\text{true positives}+\text{false positives}},\; \\
		\end{equation*}
		\begin{equation*}
			\text{recall} = \frac{\text{true positives}}{\text{true positives}+\text{false negatives}}.
		\end{equation*}
	}
\end{itemize}

\subsection{Alzheimer's Disease diagnosis}
\label{sec:5.1}

The Open Access Series of Imaging Studies (OASIS) is a project aimed at making neuroimaging data sets of the brain freely available for the scientific community. In particular, OASIS-3 is a compilation of MRI and PET imaging and related clinical data for 1098 participants who were collected across several ongoing studies in the Washington University Knight Alzheimer Disease Research Center over the course of 15 years. Imaging data is accompanied by dementia and APOE status and longitudinal clinical and cognitive outcomes \cite{lamontagne2019oasis}.

We consider a subset of the full study group, in order to have a balanced set of data.

A summary of demographic and neuropsychological details of the subjects considered in our study is presented in Table \ref{tab:subj_data}. 

\begin{table}[H]
	\centering
	\begin{tabularx}{\textwidth}{@{}l *4{>{\centering\arraybackslash}X}@{}}
		\toprule
		\quad & AD (mean) & AD (st.dev.) & Control (mean) & Control (st.dev.) \\
		\midrule
		\small{No. of subjects} & 225 & -  & 248 & - \\
		\small{Gender (F/M)} & 114/111 & -  & 126/122 & - \\
		\small{Hand preference (A/L/R)} & 5/23/197 & -  & 6/26/216 & - \\
		\small{Age at entry} & 74.41 & 7.60  & 65.21 & 9.62\\
		\small{Education (years)} & 14.77 & 3.08  & 16.04 & 2.51\\
		\small{MMSE} & 20.33 & 6.38  & 29.27 & 1.30\\
		\bottomrule
	\end{tabularx}
	\caption{Demographic details and baseline cognitive status measures of the study population.}
	\label{tab:subj_data}
\end{table}

For each subject, we build the persistence diagrams using the estimation of cortical thickness on $34$ points in both right and left hemisphere of the brain, for a total of $64$ values. For simplicity, in the study we consider the same coordinates of the above mentioned points for all subjects. The coordinates are computed with the \textit{scipy} toolbox \cite{2020SciPy-NMeth}. From this coordinates we build the persistence diagrams and we extract $1$ and $2$-dimensional topological features, i.e., we obtain the generators associated with $H_1$ and $H_2$ homological groups.

\begin{figure}[H]
	\centering
	\includegraphics[width=.48\textwidth]{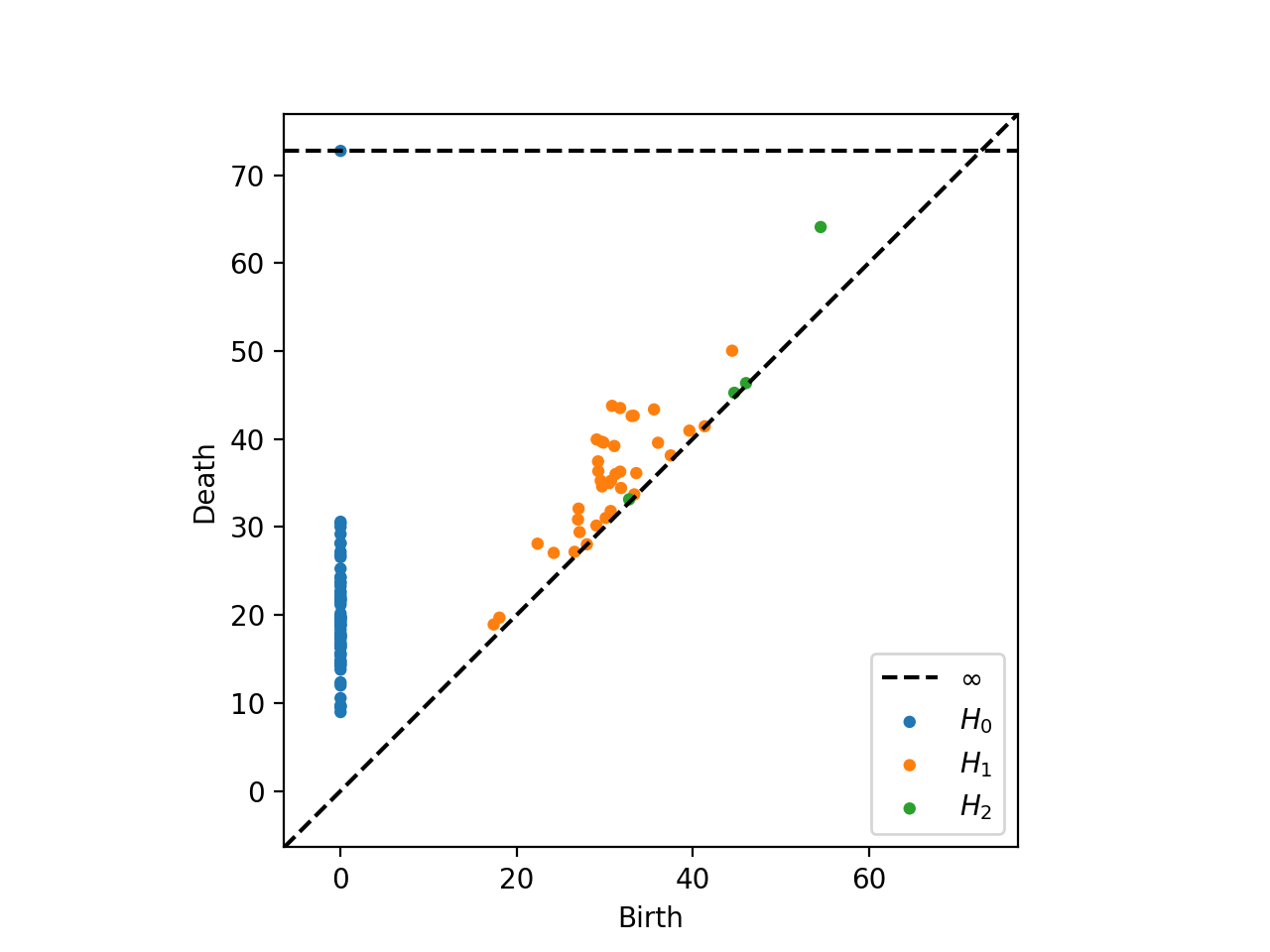} \,
	\includegraphics[width=.48\textwidth]{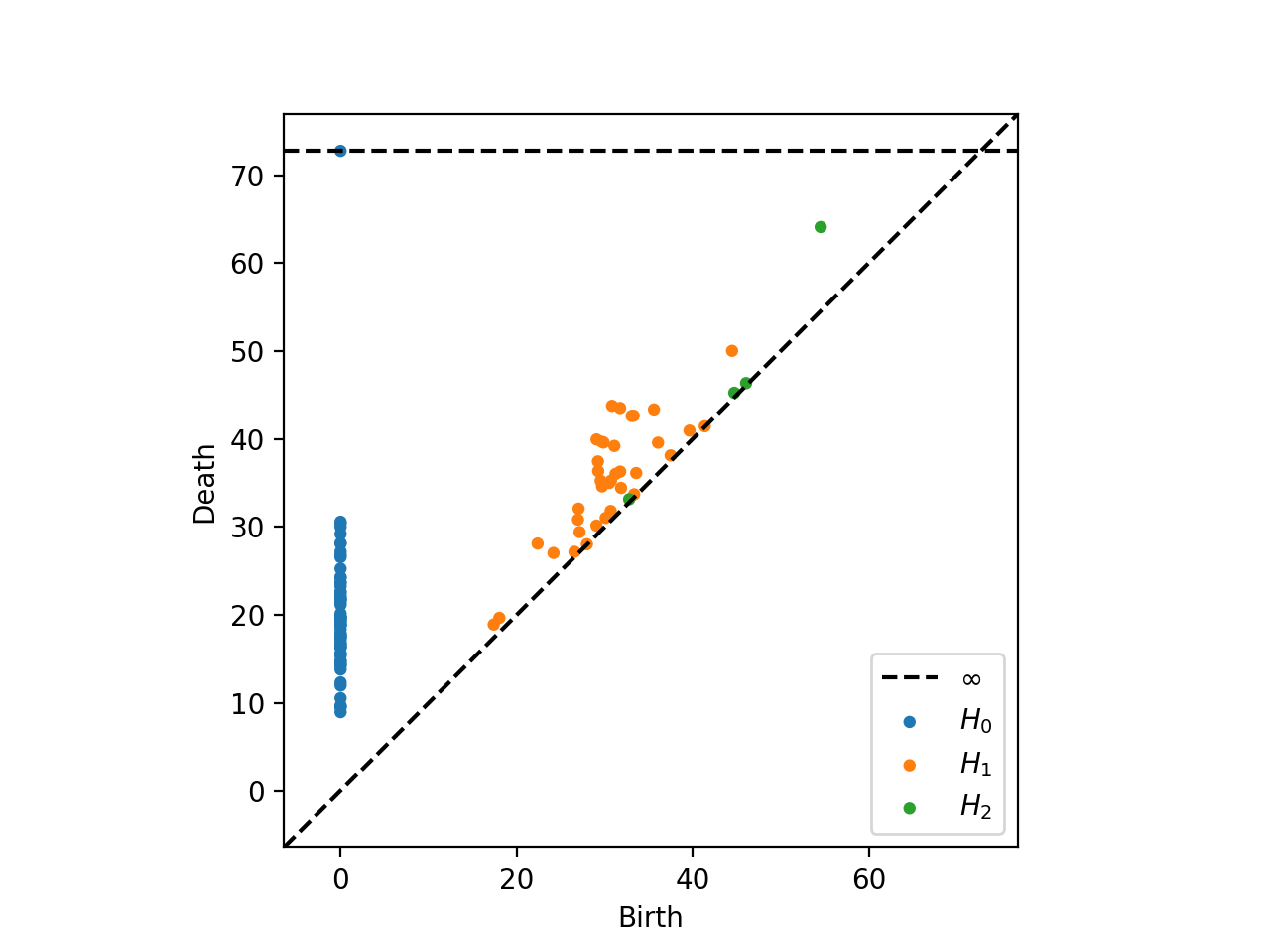}
	\caption{Persistence diagrams of an AD subject with a MMSE of $30$ (left) and the persistence diagram of a control subject with a MMSE of $7$ (right).}
	\label{fig:pd_adc}
\end{figure}

\begin{figure}[H]
	\centering
	\includegraphics[width=.48\textwidth]{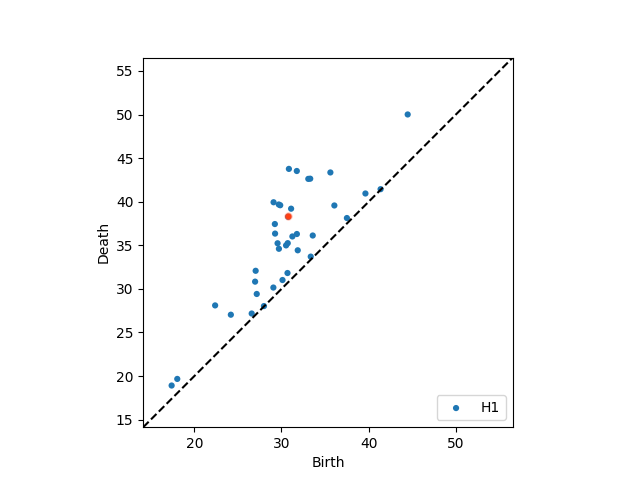}\,
	\includegraphics[width=.48\textwidth]{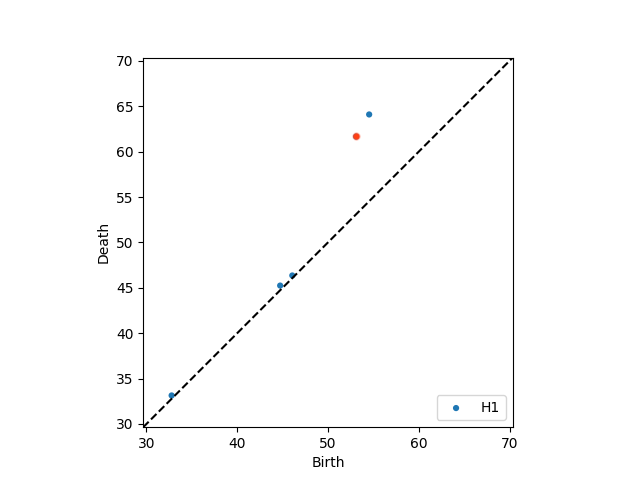}
	\caption{$1$-dimensional (left) and $2$-dimensional (right) persistence diagram of an AD subject. The red dot is the added centre of persistence via $Psi_a$ and $\psi_2$.}
	\label{fig:pd_v}
\end{figure}

In Figure \ref{fig:pd_adc} we show two examples of persistence diagrams, and in Figure \ref{fig:pd_v} we highlight the generator added as centre of persistence.

We evaluate the performance achieved by a SVMs classifier that makes use of the presented PSS, PWG and SW kernels, both in the classical and in the variably scaled settings.

In each test, we perform a random $70\%/30\%$ splitting of the dataset for training and testing, and we consider $5$-fold cross validation on the training set for the tuning of the hyperparameters. The results displayed in Tables \ref{tab:ker_res} \ref{tab:ker_w_res} \ref{tab:ker_s_res} have been averaged over $10$ runs of tests. 

For the variably scaled setting, we consider $\Psi_a$ for both $H_1$ and $H_2$ diagrams, while $\Psi_{\rho}$, with $\rho=10$, is employed with $H_1$ only, since $H_2$ diagrams are limited in the number of generators, and therefore compressing features is not meaningful.

Furthermore, we use $\psi_2$ as auxiliary function. Indeed, as highlighted in Table \ref{tab:cop}, $\psi_2$ definitely outperforms $\psi_1$ in our setting. Moreover, we observe that the performances achieved by the auxiliary function alone, i.e., taking the centres of mass or persistence in place of the persistence diagram, are definitely not competitive with respect to the classical and variably scaled settings.

\begin{table}[H]
	\centering
	\begin{tabular}{lccccc}
		\toprule
		\quad & $\Psi$ & $\psi$ & Accuracy & f$_1$-score & validation time (s)\\
		\midrule
		SW ($H_1$) & - & - & $0.741$    & $0.716$    & $220$         \\
		VSP-SW ($H_1$) & $\Psi_a$ & $\psi_2$ & $0.732$    & $0.700$      & $223$      \\
		VSP-SW ($H_1$) & $\Psi_{\rho}$ & $\psi_2$ & $0.731$   & $0.693$   & $184$     \\
		SW ($H_2$)       & - & - & $0.741$    & $ 0.712$      & $183$         \\
		VSP-SW ($H_2$)    & $\Psi_a$ & $\psi_2$ & $  0.753$    & $ 0.720$      & $300$         \\
		\bottomrule
	\end{tabular}
	\caption{OASIS-3 dataset. Results of SVMs classification obtained considering $H_1$ and $H_2$ persistence diagrams and using the SW kernel.}
	\label{tab:ker_res}
\end{table}

\begin{table}[H]
	\centering
	\begin{tabular}{lccccc}
		\toprule
		\quad & $\Psi$ & $\psi$ & Accuracy & f$_1$-score & validation time (s)\\
		\midrule
		PWGK ($H_1$) & - & - & $0.749$    & $0.723$    & $18165$         \\
		VSP-PWG ($H_1$) & $\Psi_a$ & $\psi_2$  & $0.750$    & $0.726$   & $18082$         \\
		VSP-PWG ($H_1$) & $\Psi_{\rho}$ & $\psi_2$ & $0.759$   & $0.735$   & $8731$   \\
		PWG ($H_2$) & - & -  & $0.716$    & $ 0.699$   & $4422$         \\
		VSP-PWG ($H_2$) & $\Psi_a$ & $\psi_2$   & $  0.709$    & $ 0.683$    & $4713$         \\
		\bottomrule
	\end{tabular}
	\caption{OASIS-3 dataset. Results of SVMs classification obtained considering $H_1$ and $H_2$ persistence diagrams and using the PWG kernel.}
	\label{tab:ker_w_res}
\end{table}

\begin{table}[H]
	\centering
	\begin{tabular}{lccccc}
		\toprule
		\quad & $\Psi$ & $\psi$ & Accuracy & f$_1$-score & validation time (s)\\
		\midrule
		PSS ($H_1$)& - & - & $0.743$    & $0.721$    &$9238$         \\
		VSP-PSS ($H_1$) & $\Psi_a$ & $\psi_2$  & $0.752$    & $0.728$    & $9045$         \\
		VSP-PSS ($H_1$) & $\Psi_{\rho}$ & $\psi_2$ & $0.750$   & $0.723$   & $3330$    \\
		PSS ($H_2$) & - & - & $0.781$    & $ 0.762$    &$2825$         \\
		VSP-PSS ($H_2$) & $\Psi_a$ & $\psi_2$   & $  0.775$    & $ 0.755$    &$3120$         \\
		\bottomrule
	\end{tabular}
	\caption{OASIS-3 dataset. Results of SVMs classification obtained considering $H_1$ and $H_2$ persistence diagrams and using the PSS kernel.}
	\label{tab:ker_s_res}
\end{table}

\begin{table}[H]
	\centering
	\begin{tabular}{lcc}
		\toprule
		\quad                   & accuracy    &  f$_1$-score\\ 
		\midrule
		$\psi_1$                &  $0.56$   &  $0.55$     \\
		$\psi_2$             &  $0.65$  &  $0.72$     \\
		\bottomrule
	\end{tabular}
	\caption{OASIS-3 dataset. Results of SVMs classification obtained by using the centre of mass and persistence alone in place of the persistence diagrams.}
	\label{tab:cop}
\end{table}

We observe that VSPKs are competitive with respect to the classical setting, improving the performance in some cases. Moreover, the usage of $\Psi_{\rho}$ leads to a consistent saving in validation time.

\subsection{Orbit recognition}

\begin{figure}[H]
	\centering
	\includegraphics[width=.32\textwidth]{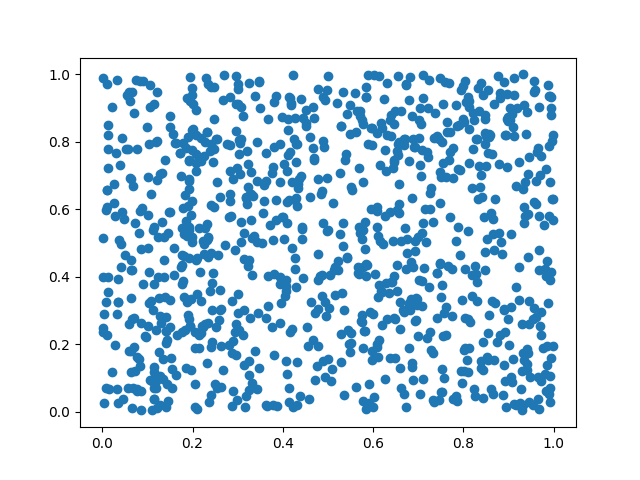}\,
	\includegraphics[width=.32\textwidth]{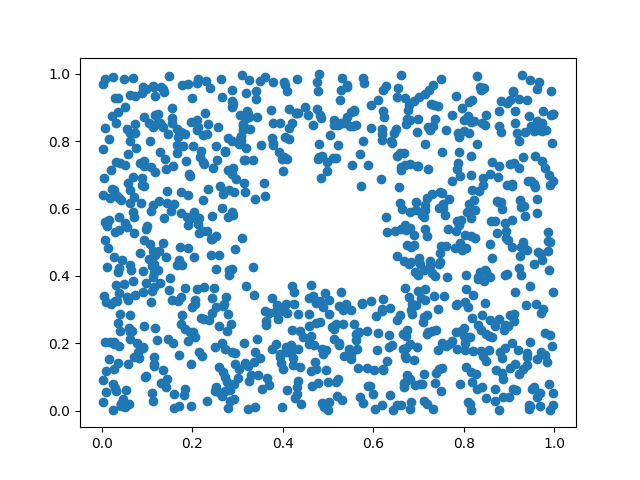}\,
	\includegraphics[width=.32\textwidth]{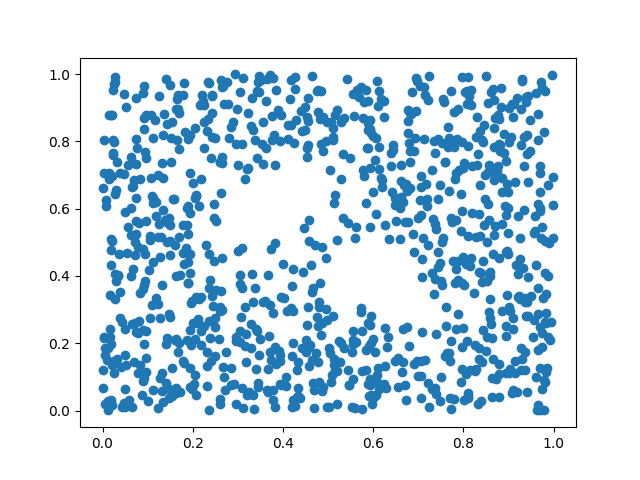}
	\caption{Fixed $(x_0, y_0) \in [0,1]^2$, the orbits resulting from the linked twisted map taking $r = 2.5, 4.1, 4.3$, from left to right, respectively.}
	\label{fig:orbit}
\end{figure}

\begin{figure}[H]
	\centering
	\includegraphics[width=.32\textwidth]{Figure_4.3.png}\,
	\includegraphics[width=.32\textwidth]{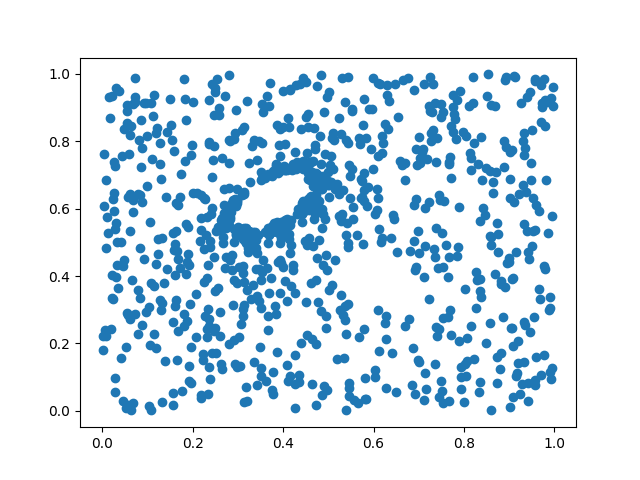}\,
	\includegraphics[width=.32\textwidth]{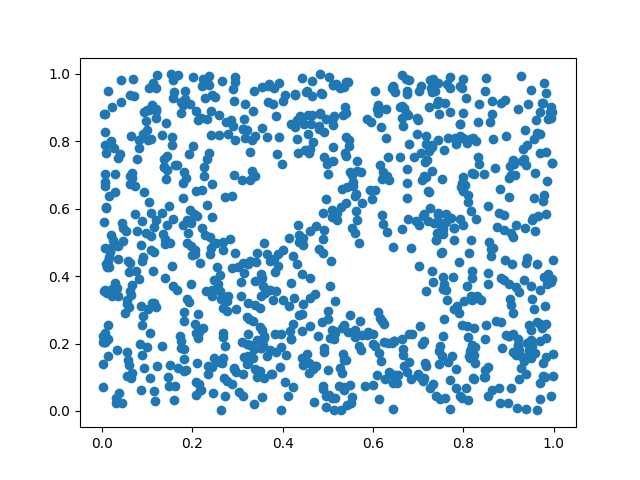}
	\caption{Fixed $r = 4.3$, the orbits resulting from the linked twisted map taking different starting points $(x_0, y_0) \in [0,1]^2$.}
	\label{fig:orbit_fixed_r}
\end{figure}

As second experiment, we follow the idea proposed in \cite{adams2017persistence} and we analyse the linked twisted map, which models fluid flows. The corresponding orbits are computed via the discrete system
\begin{equation*}
	\begin{cases}
		x_{n+1} = x_n + r y_n (1-y_n) \quad \mod 1\\
		y_{n+1} = y_n + r x_{n+1} (1-x_{n+1}) \quad \mod 1
	\end{cases}
\end{equation*}
where $(x_0, y_0) \in [0,1]\times[0,1]$ is the initial position and $r > 0$ is a real parameter that influences the orbit. The topological structure of the orbit changes with the initial position and $r$, as displayed in Figures \ref{fig:orbit} and \ref{fig:orbit_fixed_r}, where we depict the first $1000$ \textit{iterations} $\{(x_n, y_n) \, : \, n = 0, \dots, 1000\}$.

In the following tests, accordingly to \cite{adams2017persistence}, we choose a set of five parameters $r = 2.5, 3.5, 4, 4.1, 4.3$ as set of classification labels. For each label, we compute the first $1000$ points of $50$ orbits, with random starting point. Therefore, the dataset consists of $250$ elements. Then, we compute the persistence diagram related to each orbit.

Here, since each persistence diagram has a huge number of generators ($\approx 10^5$), we empirically reduce the persistence diagrams by restricting to the most $10$ persistent elements. More precisely, here $\Psi_a$ is computed with respect to such $10$ elements, while in the case of $\Psi_{\rho}$ the discarded less persistent generators are compressed in a unique element via $\psi_2$. Moreover, since we are dealing with $2$-dimensional orbits, we compute the $H_1$ homology group only.

As in Section \ref{sec:5.1}, we consider a $5$-fold cross validation on the training set, and the results displayed in Tables \ref{tab:ker_dyn}, \ref{tab:ker_w_dyn} and \ref{tab:ker_s_dyn} are averaged over $10$ runs with 70\%/30\% training-test split of the data.

\begin{table}[H]
	\centering
	\begin{tabular}{lcccc}
		\toprule
		\quad & $\Psi$ & $\psi$ & Accuracy & f$_1$-score\\
		\midrule
		SW ($H_1$) & - & -  & $0.832$    & $0.831$\\
		VSP-SW ($H_1$) & $\Psi_a$ & $\psi_2$  & $0.814$    & $0.812$\\
		VSP-SW ($H_1$) & $\Psi_{\rho}$ & $\psi_2$  & $0.833$    & $0.832$\\
		\bottomrule
	\end{tabular}
	\caption{Orbit Recognition. Results of SVMs classification on $H_1$ persistence diagrams using the SW kernel.}
	\label{tab:ker_dyn}
\end{table}

\begin{table}[H]
	\centering
	\begin{tabular}{lcccc}
		\toprule
		\quad & $\Psi$ & $\psi$ & Accuracy & f$_1$-score\\
		\midrule
		PWG ($H_1$) & - & -  & $0.858$    & $0.866$\\
		VSP-PWG ($H_1$) & $\Psi_a$ & $\psi_2$   & $0.858$    & $0.866$ \\
		VSP-PWG ($H_1$) & $\Psi_{\rho}$ & $\psi_2$ & $0.846$    & $0.853$  \\
		\bottomrule
	\end{tabular}
	\caption{Orbit Recognition. Results of SVMs classification on $H_1$ persistence diagrams using the PWG kernel.}
	\label{tab:ker_w_dyn}
\end{table}

\begin{table}[H]
	\centering
	\begin{tabular}{lcccc}
		\toprule
		\quad & $\Psi$ & $\psi$ & Accuracy & f$_1$-score\\
		\midrule
		PSS ($H_1$)  & - & -  & $0.806$    & $0.803$\\
		VSP-PSS ($H_1$)  & $\Psi_a$ & $\psi_2$  & $0.824$    & $0.821$\\
		VSP-PSS ($H_1$) & $\Psi_{\rho}$ & $\psi_2$ & $0.826$    & $0.823$\\
		\bottomrule
	\end{tabular}
	\caption{Orbit Recognition. Results of SVMs classification on $H_1$ persistence diagrams using the PSS kernel.}
	\label{tab:ker_s_dyn}
\end{table}

\subsection{3D shape segmentation}

\begin{figure}[H]
	\centering
	\includegraphics[width=.3\textwidth]{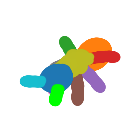}\,
	\includegraphics[width=.3\textwidth]{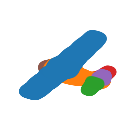}\\
	\includegraphics[width=.3\textwidth]{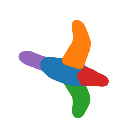}\,
	\includegraphics[width=.3\textwidth]{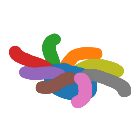}\\
	\caption{Examples of 3D shapes and corresponding segmentation labels. From top left, clock-wise, Ant, Airplane, Bird and Octopus.}
	\label{fig:seg}
\end{figure}

Here, we follow an experiment proposed also in \cite{carriere2017sliced}. We consider some categories of the mesh segmentation benchmark introduced in \cite{Chen:2009:ABF}, which contains different 3D shapes of several categories. In each sample, every face is represented as a triplet of 3D points, and it is labeled with respect to a corresponding segmentation index. In our dataset, we consider $72546$ faces in the category \textit{Ant}, $59546$ faces in the category \textit{Airplane}, $50902$ faces in the category \textit{Bird} and $66248$ faces in the category \textit{Octopus}. A $H_1$ persistence diagram is then computed for each face by using the geodesic distance on the 3D shape; for more details concerning the dataset, we refer to \cite{carriere2015stable}, while in Figure \ref{fig:seg} we display some examples of shape segmentation. 

For each category, the classification task consists in predicting the segmentation index corresponding to the persistence diagram. As in the previous subsections, we consider a $5$-fold cross validation on the training set, and the results displayed in Tables \ref{tab:seg_acc} and \ref{tab:seg_f1} are averaged over $10$ runs with 70\%/30\% training-test split of the data.

We remark that the elements of the persistence diagrams in this experiment are limited in numbers and free of noisy generators, because of the particular construction of the diagrams obtained via the geodesic metric. Therefore, the usage of $\Psi_{\rho}$ is not significant.

\begin{table}[H]
	\centering
	\begin{tabular}{lcccc}
		\toprule
		\quad & PSS ($H_1$)  & VSP-PSS ($H_1$)  & SW ($H_1$)  & VSP-SW ($H_1$) \\
		\midrule
		Ant          & $0.781$    & $0.793$    & $0.753$    & $0.746$   \\
		Airplane     & $0.677$    & $0.688$    & $0.717$    & $0.704$   \\
		Bird         & $0.612$    & $0.598$    & $0.636$    & $0.629$   \\
		Octopus      & $0.771$    & $0.773$    & $0.755$    & $0.746$   \\
		\bottomrule
	\end{tabular}
	\caption{Accuracy achieved by the SVMs classifiers in the carried out tests. In the variably scaled setting, we used $\Psi_a$ and $\psi_2$ functions.}
	\label{tab:seg_acc}
\end{table}

\begin{table}[H]
	\centering
	\begin{tabular}{lcccc}
		\toprule
		\quad & PSS ($H_1$)  & VSP-PSS ($H_1$)  & SW ($H_1$)  & VSP-SW ($H_1$) \\
		\midrule
		Ant          & $0.744$    & $0.751$    & $0.710$    & $0.700$   \\
		Airplane     & $0.580$    & $0.582$    & $0.580$    & $0.563$   \\
		Bird         & $0.619$    & $0.616$    & $0.624$    & $0.617$   \\
		Octopus      & $0.704$    & $0.707$    & $0.679$    & $0.666$   \\
		\bottomrule
	\end{tabular}
	\caption{f1-score achieved by the SVMs classifiers in the carried out tests. In the variably scaled setting, we used $\Psi_a$ and $\psi_2$ functions.}
	\label{tab:seg_f1}
\end{table}

\section{Conclusions and future work}
\label{sec:6}
In this paper, we proposed VSPKs for dealing with persistence diagrams in the context of persistent homology. The proposed framework, which is directly inspired by the variably scaled setting explored in kernel-based approximation and learning, may enhance the performance and the efficiency of existing kernels for persistence diagrams, as suggested by the obtained results. Future work consists of investigating more on the design of the scaling function, which plays a key role in the construction of the kernel. In this view, the analysis of VSPKs in the context of algebraic topology may provide useful insights.

\section*{Acknowledgements}
This research has been accomplished within the Rete ITaliana di Approssimazione (RITA) and the thematic group on Approximation Theory and Applications of the Italian Mathematical Union. We also received the support of GNCS-IN$\delta$AM.
	
\bibliographystyle{elsarticle-num}
\bibliography{bibliography.bib}
	
\end{document}